\documentstyle[12pt]{article}
\title{Chemical Examples in Hypergroups}
\author{ {\bf B. Davvaz and A. Dehghan-Nezhad } \\
Department of Mathematics, University of Yazd,\\
Yazd, Iran \\
E-mail: {\it davvaz@yazduni.ac.ir}}
\date{}
\begin{document}
\maketitle
\begin{abstract}
%
%
Hypergroups first were introduced by Marty in 1934. Up to now many
researchers have been working on this field of modern algebra and developed it.
It is purpose of this paper to provide examples of hypergroups associated
with chemistry. The examples presented are connected to construction from chain
reactions.
%
%
\end{abstract}
%
\section{Introuduction}
The theory of algebraic hyperstructures which is a generalization of
the concept of algebraic structures first was introduced by Marty in
1934 [4], and had been studied in the following decades and nowadays
by many mathematicians, and many papers concerning various
hyperstructures have appeared in the literature, for example see
[2,3,6,8]. The basic definitions of the object can be found in
[1,7].\\
\subsection{Hypergroups and $H_v$-groups}

 An {\it algebraic  hyperstructure} is a non-empty set $H$ together with a function
$\cdot : H \times H \longrightarrow {p}^*(H)$ called {\it
hyperoperation}, where ${p^*(H})$ denotes the set of all non-empty
subsets of $H$. If $A,B$ are non-empty subsets of $H$ and $x \in H$,
then we define
$$ A\cdot B=\displaystyle \bigcup_{a\in A,b\in B}
a\cdot b,\ \ x \cdot B=\{ x \} \cdot B, \ \ {\rm and} \ \ A \cdot
x=A \cdot \{ x \}. $$

The hyperoperation $(\cdot )$ is called {\it associative} in $H$ if
$$
(x \cdot y) \cdot z =  x \cdot (y \cdot z)  \ {\rm for \ all} \
x,y,z \ {\rm in} \ H, $$
 which means that $$\displaystyle
\bigcup_{u\in x\cdot y}u\cdot z=\bigcup_{v\in y\cdot z}x\cdot v.$$

We say that a  semihypergroup $(H,\cdot)$ is a {\it hypergroup} if
for all $x\in H$, we have $x\cdot H=H\cdot x=H.$ A hypergroupoid
$(H,\cdot )$ is an {\it $H_v$-group}, if for all $x,y,z\in H,$ the
following conditions hold:
\begin{itemize}
\item [(1)] $x\cdot (y\cdot z) \cap (x\cdot y)\cdot z \not =
\emptyset$, (weak associative) \item [(2)] $x\cdot H=H\cdot x=H$.
\end{itemize}

A non-empty subset $K$ of a hypergroup (respectively, $H_v$-group)
$H$ is called a {\it subhypergroup} (respectively, $H_v$-subgroup)
of $H$ if $a \cdot K=K \cdot a=K$ for all $a \in K$.
\\
\indent In this paper, we will give some examples of hypergroups associated with chemistry.
The examples presented are connected to construction from chian reactions.
%
\section{Prelimiaries}
{\bf a) Chain reactions}\\
\indent An atom of group of atoms possessing an odd (unpaired) electron is called a
free radical, such as
\[
Cl, \ CH_3, \ C_2 H_5
\]
The chlorination of methane is an example of a chain reaction, a reaction that
involves a series of setps, each of which generates a reactive substance that
brings about the next step. While chain reactions may vary widely in their
details, they all have certain fundamental characteristics in common.
\begin{itemize}
\item[1)]$Cl_2 \longrightarrow 2Cl^o $ \\
(1) is called Chain-initiating step.
\item[2)]$Cl^o+CH_4 \longrightarrow HCl+CH_3^o $
\item[3)]$CH_3^o+Cl_2 \longrightarrow CH_3Cl+Cl^o $   \\
then (2), (3), (2), (3), etc, until finatly: \\
(2) and (3) are called Chain-propagating steps.
\item[4)]$Cl^o+Cl^o \longrightarrow Cl_2 \ \ \ $ or
\item[5)]$CH_3^o+CH_3^o \longrightarrow CH_3CH_3 \ \ \ $  or
\item[6)]$CH_3^o+Cl^o \longrightarrow CH_3Cl $.  \\
(4),(5) and (6) are called Chain-terminating steps.
\end{itemize}

First in the chain of reactions is a chain-initiating step,
in which energy is
absorbed and a reactive particle generated; in the present reaction
it is the cleavage of chlorine into atoms (step 1).\\

There are one or more chain-propagating steps, each of which consumes a
reactive particle and generates another; there they are the reaction of
chlorine atoms with methane (step 2), and of methyl radicals with chlorine
(step 3).\\

Finally, there are chain-terminating steps, in which reactive particles are
consumed but not generated; in the chlorination of methane these would involve
the union of two of the reactive particles, or the capture of one of them by
the walls of the reaction vessel.\\  \\
{\bf b) The Halogens F, CL, Br, and I}\\
\indent  The halogens are all typical non-metals. Although their physical forms
differ-fluorine and chlorine are gases, bromine is a liquid and iodine is a solid at room
temprature, each consists of diatomic molecules; $F_2, Cl_2, Br_2$ and $I_2$.
The halogens all react with hydrogen to form gaseous compounds, with the
formulas
$HF, HCL, HBr, $ and $HI$
all of which are very soluble in water.
The halogens all react with metals to give halides.
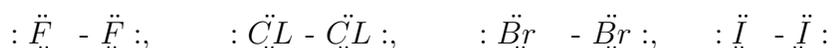
\begin{center}
$:\ddot {F}$\hspace*{-2mm}$_{..}$ \hspace*{1mm} - $\ddot {F}:,$\hspace*{-6mm}$_{..}$ \hspace*{12mm}
$:\ddot {CL}$\hspace*{-4mm}$_{..}$ \hspace*{1mm} - $\ddot {CL}:,$\hspace*{-8mm}$_{..}$ \hspace*{6mm}
$ \ \ \ \ \  $ $:\ddot {Br}$\hspace*{-3mm}$_{..}$ \hspace*{3mm} -
$\ddot {Br}:,$\hspace*{-7mm}$_{..}$ \hspace*{3mm}$ \ \ \ \ \   $
$:\ddot {I}$\hspace*{-2mm}$_{..}$ \hspace*{1mm} -
$\ddot {I}:$\hspace*{-5mm}$_{..}$ \hspace*{3mm}  \
\end{center}
\indent The reader will find in [5] a deep discussion of chain reactions and
halogens.
%
%
\section{Chemical Hypergroups}
In during chain reaction
$$
A_2+ B_2
\stackrel{\rm Heat \ or\  Light} \longleftrightarrow
2AB
$$
there exist all molecules $A_2, B_2, AB$ and whose fragment parts $A^o, B^o$ in
experiment. Elements of this colletion can by combine with each other.
All combinational probability for the set ${\cal H} = \{ A^o, B^o, A_2, B_2, AB \}$
to do without energy can be displayed as follows:
{\tiny
$$
\begin{tabular}{c|c|c|c|c|c}
$+$   & $A^o           $ & $B^o           $  & $ A_2              $ & $B_2              $ & $AB$  \\  \hline
$A^o$ & $A^o,A_2       $ & $A^o,B^o,AB    $  & $ A^o, A_2         $ & $A^o,B_2,B^o,AB   $ & $A^o,AB,A_2,B^o$ \\  \hline
$B^o$ & $A^o,B^o,AB    $ & $B^o,B_2       $  & $ A^o,B^o,AB,A_2   $ & $B^o,B_2          $ & $A^o,B^o,AB,B_2$  \\  \hline
$A_2$ & $A^o,A_2       $ & $A^o,B^o,AB,A_2$  & $A^o, A_2          $ & $A^o,B^o,A_2,B_2,AB$& $A^o,B^o,A_2,AB $  \\  \hline
$B_2$ & $A^o,B^o,B_2,AB$ & $B^o,B_2       $  & $A^o,B^o,A_2,B_2,AB$ & $B^o, B_2         $ & $A^o,B^o,B_2,AB$  \\  \hline
$AB$  & $A^o,AB,A_2,B^o$ & $A^o,B^o,AB,B_2$  & $A^o,B^o,A_2,AB    $ & $A^o,B^o,B_2,AB   $ & $A^o,B^o,A_2,B_2,AB$\\  \hline
\end{tabular}
$$
}
\normalsize
{\bf Theorem.} $({ \cal H} , + )$  is an $H_v$-group. \\ \\
{\bf Proof.} Clearly reproduction axiom and weak associativity are valid. As a sample
of how to calculate the weak associativity, we illustrate some cases:
$$
\left\{
\begin{array}{l}
(AB+A_2)+B_2 = \{ AB, A_2, A^o, B^o \}+B_2 = \{B_2, AB, A_2, A^o, B^o \},   \\
AB+(A_2+B_2)=AB+\{A_2, B_2, A^o, B^o, AB \} = \{ A_2, B_2, AB, A^o, B^o\},
\end{array}
\right.
$$
$$
\ \ \ \left\{
\begin{array}{l}
(AB+A^o)+A^o=\{ AB, A^o, A_2, B^o \}+A^o = \{A_2, A^o, AB, B^o \}, \\
AB+(A^o+A^o)=AB+ \{A_2, A^o \} = \{ A_2, AB, A^o, B^o\},
\end{array}
\right.
$$
$$
\left\{
\begin{array}{l}
(A_2+B^o)+B_2 = \{ AB, A^o, A_2, B^o \}+B_2 = \{B_2, AB, B^o, A^o, A_2 \}, \\
A_2+(B^o+B_2)=A_2+\{B_2, B^o \} = \{ A_2, A^o, AB, B^o, B_2\}.
\end{array}
\right.
$$

\noindent {\bf Corollary.} ${\cal H}_1=\{A^0,A_2 \}$ and ${\cal H}_2=\{B^0,B_2 \}$ are
only subhypergroups of $({\cal H}, +)$.\\ \\
\indent If we consider $A=H$ and $B \in \{ F, CL, Br, I \}$ (for example $B=I$),
the complete reaction table becomes:
{\tiny
$$
\begin{tabular}{c|c|c|c|c|c}
$+$   & $H^o           $ & $I^o           $  & $ H_2              $ & $I_2              $ & $HI$  \\  \hline
$H^o$ & $H^o,H_2       $ & $H^o,I^o,HI    $  & $ H^o, H_2         $ & $H^o,I_2,I^o,HI   $ & $H^o,HI,H_2,I^o$ \\  \hline
$I^o$ & $H^o,I^o,HI    $ & $I^o,I_2       $  & $ H^o,I^o,HI,H_2   $ & $I^o,I_2          $ & $H^o,I^o,HI,I_2$  \\  \hline
$H_2$ & $H^o,H_2       $ & $H^o,I^o,HI,I_2$  & $H^o, H_2          $ & $H^o,I^o,H_2,I_2,HI$& $H^o,I^o,H_2,HI $  \\  \hline
$I_2$ & $H^o,I^o,I_2,HI$ & $H^o,I_2       $  & $H^o,I^o,H_2,I_2,HI$ & $H^o, I_2         $ & $H^o,I^o,I_2,HI$  \\  \hline
$HI$  & $H^o,HI,H_2,I^o$ & $H^o,I^o,HI,I_2$  & $H^o,I^o,H_2,HI    $ & $H^o,I^o,H_2,HI   $ & $H^o,I^o,H_2,I_2,HI$\\  \hline
\end{tabular}
$$
}
\normalsize
%
%
{\bf Acknowledgment} \\
\indent We appreciate the assistance and suggestions of Dr. A. Gorgi at the
Department of Chemistry.
%
%


\newpage
\begin{thebibliography}{99}
\bibitem{1} P. Corsini, {\it Prolegomena of hypergroup theory}, second
            edition, Aviani editor, (1993).
\bibitem{2} M.R. Darafsheh and B. Davvaz, {\it $H_v$-ring of fractions},
             Italian J. Pure Appl. Math., 5 (1999) 25-34.
\bibitem{3} B. Davvaz, {\it Weak Polygroups}, Proc. 28th Annual Iranian Math.
Conf, (1977), 139-145.
\bibitem{4} F. Marty, {\it Sur une generalization de la notion de groupe}, $8^{iem}$
congres Math. Scandinaves, Stocklhom, (1934), 45-49.
\bibitem{5} Morrison and Boyd, {\it Organic Chemistry}, Sixth Eddition,
Prentice-Hall, Inc, 1992.
\bibitem{6} T. Vougiouklis, {\it A new class of hyperstructures}, J. Combin.
Inf. system Sci, 20, (1995), 229-235.
\bibitem{7} T. Vougiouklis, {\it Hyperstructures and their repesentations},
Hadronic Press Inc, (1994).
\bibitem{8} T. Vougiouklis, {\it Convolution on WASS hyperstructures},
Discrete Math., (1997), 347-355.
%
\end{thebibliography}
\end{document}